\documentclass[fleqn]{mat01}
\usepackage{times,mathtimy,amssymb}

\newtheorem{theor}{\bf Theorem}
\newtheorem{lem}{\it Lemma}
\newtheorem{propo}{\rm PROPOSITION}
\newtheorem{coro}{\rm COROLLARY}
\newtheorem{definit}{\rm DEFINITION}
\newtheorem{exa}{\it Example}
\newtheorem{rem}{\it Remark}

\newcommand{\supp}{\operatorname {Supp}}
\newcommand{\ann}{\operatorname {Ann}}
\newcommand{\ass}{\operatorname {Ass}}
\newcommand{\assh}{\operatorname {Assh}}

\newcommand{\lth}{\operatorname {length}}

\newcommand{\dth}{\operatorname {depth}}
\newcommand{\ov}{\bar}
\newcommand{\m}{\mathfrak m}
\newcommand{\V}{\operatorname {V}}
\newcommand{\cm}{Cohen--Macaulay}
\newcommand{\biss}{,\dots ,}

\renewcommand{\labelenumi}{$(\roman{enumi})$}

\begin{document}

\setcounter{page}{159} \firstpage{159}

\font\zxzx=msam10 at 10pt
\def\twoheadrightarrow{\mbox{\zxzx{\char'263\ \ }}}

\renewcommand\theequation{\thesection\arabic{equation}}

\title{Reducing system of parameters and the Cohen--Macaulay property}

\markboth{Bj\"orn M\"aurer and J\"urgen St\"uckrad}{Reducing
system of parameters}

\author{BJ\"ORN M\"AURER and J\"URGEN ST\"UCKRAD}

\address{Fakult\"at f\"ur Mathematik und Informatik, Universit\"at Leipzig, Augustus Platz 10/11, D-04109 Leipzig,
Germany\\
\noindent E-mail: bmp@gmx.de; stueckrad@math.uni-leipzig.de}

\volume{117}

\mon{May}

\parts{2}

\pubyear{2007}

\Date{MS received 20 September 2005; revised 26 June 2006}

\begin{abstract}
Let $R$ be a local ring and let ($x_1\biss x_r$) be part of a
system of parameters of a finitely generated $R$-module $M,$ where
$r < \dim_R M$. We will show that if ($y_1\biss y_r$) is part of a
reducing system of parameters of $M$ with $(y_1\biss
y_r)M=(x_1\biss x_r)M$ then $(x_1\biss x_r)$ is already reducing.
Moreover, there is such a part of a reducing system of parameters
of $M$ iff for all primes $P\in \supp M \cap V_R(x_1\biss x_r)$
with $\dim_R R/P = \dim_R M\, -r$ the localization $M_P$ of $M$ at
$P$ is an $r$-dimensional \cm\ module over $R_P$.

Furthermore, we will show that $M$ is a \cm\ module iff $y_d$ is a
non zero divisor on $M/(y_1\biss y_{d-1})M$, where $(y_1\biss
y_d)$ is a reducing system of parameters of $M$ ($d := \dim_R M$).
\end{abstract}


\maketitle

\section{Preliminaries}

In what follows, let $R$ be a local ring with maximal ideal $\m$
and let $M$ be a non zero finitely generated $R$-module of
dimension $d$. Instead of $\dim_R$, $\dth _R$, $\ass_R$,
$\supp_R$, $\ldots$ we will write $\dim, \dth, \ass, \supp,
\ldots$ for short.

We note that $\supp M/XM = \supp M \cap V(X)$, where $X$ is a
subset of $R$ and that for a prime ideal $P$ of $R$ we have $P \in
\ass M$ iff $PR_P \in \ass M_P$. Moreover we define $\assh
M:=\{P\in\ass M\mid\dim R/P=d\}$.

For undefined terminology we refer to the standard literature
(e.g. \cite{eisenbud}).

\begin{definit}\label{defrps}$\left.\right.$\vspace{.5pc}

\noindent {\rm A system of parameters $(x_1\biss x_d)$ of $M$ is
called {\it reducing}, if for all $i=1\biss d-1$ we have
\begin{equation*}
x_i\notin P \mbox{ for all } P\in\ass M/(x_1\biss x_{i-1})M \mbox{
with } \dim R/P= d-i.
\end{equation*}}
\end{definit}
\setcounter{rem}{1}
\begin{rem}\label{defredsys}
{\rm Auslander and Buchsbaum defined in \cite{auslander} a system
of parameters $(x_1\biss x_d)$ of $R$ to be a reducing system of
parameters of $M$ if
\begin{align*}
e_{M}(x_1\biss x_d) &=\lth(M/(x_1\biss x_d)M)\\[.5pc]
&\quad\,-\lth((x_1\biss x_{d-1})M:x_d/(x_1\biss x_{d-1})M).
\end{align*}
This definition is equivalent to the definition given above if we
pass from $R$ to \hbox{$\ov R:=R/\ann_RM$} and consider $M$ as a
$\ov R$-module in Definition~\ref{defrps} and use Corollary~4.8 in
\cite{auslander}. Therefore it is clear that all definitions and
results on reducing systems of parameters remain true in this more
general context.}
\end{rem}

\begin{rem}\label{rem2}
{\rm For every system of parameters $(x_1\biss x_d)$ of $M$ there
is a reducing system of parameters $(y_1\biss y_d)$ of $M$ such
that $(y_1\biss y_d)R = (x_1\biss x_d)R$, in particular,
$(y_1\biss y_d)M = (x_1\biss x_d)M$ (see Proposition 4.9 in
\cite{auslander}).}
\end{rem}

\setcounter{definit}{3}
\begin{definit}\label{defrtps}$\left.\right.$\vspace{.5pc}

\noindent {\rm A sequence $x_1\biss x_r$ of elements of $\m$ is
called {\it part of a {\rm (}reducing{\rm )} system of parameters}
of $M$, if there are elements $x_{r+1}\biss x_d\in\m$ such that
$(x_1\biss x_r,x_{r+1}\biss x_d)$ is a (reducing) system of
parameters of $M$.}
\end{definit}

\setcounter{rem}{4}
\begin{rem}\label{remark}$\left.\right.$
{\rm \begin{enumerate}\leftskip .15pc
\renewcommand{\labelenumi}{(\arabic{enumi})}
\item A sequence $(x_1\biss x_r)$ of elements of $\m$ with $r < d$ is part of a system of parameters of $M$ iff $\dim M/(x_1\biss x_r)M = d-r$.
\item A sequence $(x_1,\dots,x_r)$ of elements of $\m$ with $r < d$ is part of a reducing system of parameters of $M$ iff for all $i=1\biss r$ we have
$x_i \notin P$ for all $P\in\ass M/(x_1\biss x_{i-1})M$ with $\dim
R/P \ge d-i$.
\item Every regular sequence on $M$ is part of a reducing system of parameters of $M$.\vspace{-.5pc}
\end{enumerate}}
\end{rem}

\begin{rem}\label{4}$\left.\right.$

{\rm \begin{enumerate}\leftskip .15pc
\renewcommand{\labelenumi}{(\arabic{enumi})}
\item We note that the following conditions are equivalent:
\leftskip .5pc \vspace{6pt}

(i)$M$ is a \cm\ module, i.e. $\dth M = d$.

(ii) Every system of parameters of $M$ is a regular sequence on
$M$.

(iii)There exists a system of parameters of $M$ which is a regular
sequence on $M$.\vspace{1pc}
\item Assume that $M$ is a \cm\ module. If $(x_1\biss x_r)$ is part of a system of parameters of $M$ then $M/(x_1\biss x_r)M$ is unmixed, more precisely, $\dim R/P = d-r$ for all $P\in\ass M/(x_1\biss x_r)M$. Therefore for a sequence $(x_1\biss x_r)$ of elements of $\m$ the following conditions are equivalent:
\leftskip .5pc \vspace{6pt}

(i) $(x_1\biss x_r)$ is a regular sequence on $M$.

(ii) $(x_1\biss x_r)$ is part of a reducing system of parameters
of $M$.

(iii) $(x_1\biss x_r)$ is part of a system of parameters of $M$.
\vspace{1pc}
\end{enumerate}}
\end{rem}

Let $x_1\biss x_r \in \m$. If $(x_1\biss x_r)$ is a regular
sequence on $M$ then $(x_1\biss x_r)$ is a regular sequence on
$M_P$ as well for all primes $P \in \supp M \cap V(x_1\biss x_r)$.

\setcounter{lem}{6}
\begin{lem}\label{5}
Let $(x_1,\dots,x_r)$ be part of a {\rm (}reducing{\rm )} system
of parameters of $M$. Then $(x_1,\dots,x_r)$ is part of a {\rm
(}reducing{\rm )} system of parameters of $M_P$ for all primes $P
\in \supp M \cap V(x_1,\dots,x_r)$ with $\dim R/P + \dim M_P = d$.
\end{lem}

\begin{proof}
Let $P \in \supp M \cap V(x_1,\dots,x_r)$ with $\dim R/P + \dim
M_P = d$. An easy induction argument (induction on $r$) shows that
we can restrict ourselves to the case $r=1$ (and $\dim M_P \ge
2$).

Let $\mathfrak q \in \ass M_P$\: with\: $\dim R_P/\mathfrak q =
\dim M_P(\dim R_P/\mathfrak q \ge \dim M_P-1$). Then $\mathfrak q
= QR_P$ with $Q \in \ass M$, $Q \subseteq P$, and we obtain
\begin{align*}
\dim R/Q &\geq \dim R/P +\dim(R/Q)_P = \dim R/P + \dim
R_P/\mathfrak q \\[.2pc] &= \dim R/P + \dim M_P = d\\[.2pc]
(\!\!&\geq \dim R/P + \dim M_P -1 = d-1).
\end{align*}
Therefore $x_1 \notin Q$ by our assumption. But then $x_1 \notin
\mathfrak q$, i.e. $(x_1)$ is part of a (reducing) system of
parameters of $M_P$.\hfill $\Box$
\end{proof}

\begin{lem}\label{6}
If $x\in R$ is a zero divisor on $M${\rm ,} then $P\in\ass M/xM$
for all minimal primes $P\in \ass M\cap\V(x)$.
\end{lem}

\begin{proof}
Let $P\in \ass M\cap V(x)$ be minimal. Since $P \in \ass M/xM$ iff
$PR_P \in \ass M_P/xM_P$ we may assume by localizing at $P$ that
$P = \m$. Then $x\notin Q$ for all $Q\in\ass M\backslash\{ \m\}$.
Since $R$ is noetherian there is an $i\in\mathbb N^+$ such that
$0:_Mx^i=0:_Mx^j$ for all $j\geq i$. Let $U:=0:_Mx^i$. Then $U\neq
0$ (otherwise $x$ would be a non zero divisor on $M$,
contradicting our assumption).

Let $Q\in\supp M\backslash\{ \m\}$. Since $\ass M_Q =
\{Q'R_Q|Q'\in \ass M, Q' \subseteq Q\}$, we have $x \notin
\mathfrak q$ for all $\mathfrak q \in \ass M_Q$. Therefore
$U_Q=0:_{M_Q}x^i =0$ for all $Q\in\supp M\backslash\{ \m\}$, i.e.
$\supp U = \{\m\}$. Moreover,
\begin{equation*}
U:_Mx = 0:_Mx^{i+1} = 0:_Mx^i=U.
\end{equation*}
Let $\varphi\hbox{\rm :}\ U\rightarrow M/xM$ be the inclusion
$U\subseteq M$ followed by the canonical epimorphism
$M\twoheadrightarrow M/xM$. Since
\begin{equation*}
\ker\varphi =U\cap xM=x(U:_Mx)=xU, \end{equation*} $\varphi$
induces a monomorphism $U/xU\rightarrow M/xM$. Now $U/xU\neq 0$ by
Nakayama's lemma. Therefore $\emptyset \neq\ass U/xU\subseteq\ass
U=\{ \m\}$, that means $\ass U/xU=\{ \m\}$. This gives us the
existence of a monomorphism $R/\m\rightarrow U/xU\rightarrow
M/xM$. Thus $\m\in\ass M/xM$. \hfill $\Box$
\end{proof}

\begin{lem}\label{8}
Let $Q\in\supp M$ and assume that there is an $x\in\m$ with
$x\notin Q$. Then there is a $P\in\supp M$ such that $x\in P${\rm
,} $Q\subset P$ and $\dim R/P=\dim R/Q-1$.
\end{lem}

\begin{proof}
Since $(x)$ is part of a system of parameters of $R/Q$, there is a
$P \in \supp (R/Q)/$ $x(R/Q) = \supp R/(Q + xR) = V(Q+xR) \subset
V(Q)$ with $\dim R/P = \dim R/(Q+xR) = \dim R/Q - 1$. Since $0 \ne
M_Q \cong (M_P)_{QR_P}$ we have $M_P \ne 0$, i.e. $P \in \supp
M$.

\hfill $\Box$
\end{proof}

\setcounter{coro}{9}
\begin{coro}\label{cor}$\left.\right.$\vspace{.5pc}

\noindent Let $Q\in\supp M$ and let $x_1\biss x_r\in\m$. Then
there is a $P\in\supp M\cap V(x_1\biss x_r)$ such that $Q\subseteq
P$ and $\dim R/P\ge\dim R/Q\,-r$.
\end{coro}

The proof follows immediately from Lemma \ref{8} by induction on $r$.

\section{Main results}

\setcounter{theor}{10}
\begin{theor}[\!]\label{7}
Let $(y_1\biss y_d)$ be a reducing system of parameters of $M$.
$M$ is a \cm\ module iff $y_d$ is a non zero divisor on
$M/(y_1\biss y_{d-1})M$.
\end{theor}

\begin{proof}
The implication `$\Rightarrow$' is clear, since every system of
parameters in a \cm\ module is a regular sequence (see Remark
\ref{4}(1)).

We will prove the opposite implication by induction on $d$, where
the case $d=1$ is clear. Let $d \ge 2$ and assume that the
statement is true for modules with a dimension strictly less than
$d$.

Assume that $y_d$ is a non zero divisor on $M/(y_1\biss
y_{d-1})M$. By our induction hypothesis, $M/y_1M$ is a \cm\ module
and therefore it remains to show that $y_1$ is a non zero divisor
on $M$. Suppose this is not the case. Let $P$ be minimal in $\ass
M \cap V(y_1)$. Since $(y_1)$ is part of a reducing system of
parameters of $M$, we have $\dim R/P\leq d-2$. By Lemma \ref{6},
$P\in\ass M/y_1M$ and therefore $\dim R/P =\dim M/y_1M = d-1$ (see
Remark \ref{4}(2)), a contradiction.\hfill $\Box$
\end{proof}

\setcounter{lem}{11}
\begin{lem}\label{9}
Let $(x)$ be part of a system of parameters of $M$. If $d\ge
2${\rm ,} the following conditions are equivalent{\rm :}
\begin{enumerate}
\renewcommand{\labelenumi}{\rm (\roman{enumi})}
\leftskip .4pc
\item $(x)$ is part of a reducing system of parameters of $M$.
\item $M_P$ is a one-dimensional \cm\ module over $R_P$ for all $P\in\supp M \cap\V(x)$ satisfying $\dim R/P=d-1$.
\item There is a $y \in \m$ such that $(y)$ is part of a reducing system of parameters of $M$ and $y M=x M$.
\item There is a $y \in \m$ such that $(y)$ is part of a reducing system of parameters of $M$ and $\supp M \cap V(x) \subseteq
V(y)$.\vspace{-.8pc}
\end{enumerate}
\end{lem}

\begin{proof}
The implications (i) $\Rightarrow$ (iii) and (iii) $\Rightarrow$ (iv) are obvious.

(iv) $\Rightarrow$ (ii): Let $(y)$ be part of a reducing system of
parameters of $M$ with $\supp M\cap V(x) \subseteq V(y)$ and let
$P\in\supp M \cap V(x)$ with $\dim R/P=d-1$. Then $y \in P$. Now
$y \notin Q$ for all $Q \in \ass M$ with $\dim R/Q \ge d-1$ by our
assumption. Thus $P \notin \ass M$ and therefore $PR_P \notin \ass
M_P\, (\ne \emptyset)$ from which
\begin{equation*}
0 < \dth M_P \le \dim M_P\le \dim M - \dim R/P = 1,
\end{equation*}
i.e. $\dth M_P =\dim M_P = 1$.

(ii) $\Rightarrow$ (i): Let $P\in\ass M$ with $\dim R/P\geq d-1$.
If $\dim R/P=d$, then $x\notin P$ since $(x)$ is part of a system
of parameters of $M$. Let $\dim R/P=d-1$. If $x\in P$, then $M_P$
is a \cm\ module over $R_P$ with $\dim M_P = 1$. Therefore $PR_P
\notin \ass M_P$ contradicting $P \in \ass M$.

Thus $x\notin P$ for all $P\in\ass M$ with $\dim R/P\geq d-1$,
i.e. $(x)$ is part of a reducing system of parameters of $M$ by
Remark \ref{remark}(2).\hfill $\Box$
\end{proof}

\setcounter{rem}{12}
\begin{rem}\label{rem3}
{\rm Let $x_1\biss x_r,\ y_1\biss y_r$ be elements of $\m$ with
\begin{equation*}
\supp M \cap V(x_1\biss x_r) \subseteq V(y_1\biss y_r)
\end{equation*}
(which is equivalent to $\supp M/(x_1\biss x_r)M \subseteq \supp
M/(y_1\biss y_r)M$).
\begin{enumerate}
\renewcommand{\labelenumi}{(\alph{enumi})}
\leftskip .15pc
\item If $(y_1\biss y_r)$ is part of a system of parameters of $M$ then the same is true for $(x_1\biss x_r)$. This follows immediately from Remark \ref{remark}(1).
\item If $(y_1\biss y_r)$ is a regular sequence on $M$ then the same is true for $(x_1\biss x_r)$. This follows from Corollary
2 of \cite{patil}.
\end{enumerate}
The equivalence (i) $\Leftrightarrow$ (iv) of our next theorem
shows that a similar statement holds for parts of reducing systems
of parameters of $M$, provided $r < d$. (For $r=d$ this is not
true in general, see Remark \ref{rem2}.)}
\end{rem}

\setcounter{theor}{13}
\begin{theor}[\!]\label{10}
Let $(x_1\biss x_r)$ be part of a system of parameters of $M${\rm
,} where $0 \le r <d$. Then the following conditions are
equivalent{\rm :}
\begin{enumerate}
\renewcommand{\labelenumi}{\rm (\roman{enumi})}
\leftskip .45pc
\item $(x_1\biss x_r)$ is part of a reducing system of parameters of $M$.
\item $M_P$ is an $r$-dimensional \cm\ module over $R_P$ for all $P\in\supp M \cap V(x_1\biss x_r)$ satisfying $\dim R/P=\dim M\,-r$.
\item There is a part $(y_1\biss y_r)$ of a reducing system of parameters of $M$ such that $(y_1\biss y_r)M=(x_1\biss x_r)M$.
\item There is a part $(y_1\biss y_r)$ of a reducing system of parameters of $M$ such that $\supp M\cap \V(x_1\biss x_r)\subseteq V(y_1\biss
y_r)$.\vspace{-.5pc}
\end{enumerate}
\end{theor}
\begin{proof}
We use induction on $r$. For $r=0$, there is nothing to show and
for $r=1$ the statement follows from Lemma \ref{9}. So let $r\ge
2$.

The implications (i) $\Rightarrow$ (iii) and (iii) $\Rightarrow$ (iv) are obvious.

(iv) $\Rightarrow$ (ii): Let $\ov M:=M/y_1M$. Take $P\in\supp M
\cap V(x_1\biss x_r) \subseteq \supp M\cap V(y_1\biss y_r)$ with
$\dim R/P=d-r$. Since $\dim \ov M = d-1$, $\ov M_P \cong
M_P/y_1M_P$ is an ($r-1$)-dimensional \cm\ module (over $R_P$) by
the induction hypothesis ((i) $\Rightarrow$ (ii)). Therefore it is
sufficient to show that $y_1$ is a non zero divisor on $M_P$.

Suppose this is not the case. Then by Lemma \ref{6} there is a
$\mathfrak q \in \ass M_P$ with $\mathfrak q \in \ass \ov M_P$.
Therefore $\dim R_P/\mathfrak q = r-1$ by Remark \ref{4}(2). Now
$\mathfrak q = QR_P$ with $Q \in \supp \ov M = \supp M \cap
V(y_1)$ and $Q \subseteq P$. Then $Q \in \ass M$ and we have
\begin{align*}
\dim R/Q \ge \dim R/P + \dim (R/Q)_P = \dim R/P + \dim
R_P/\mathfrak q = d-1.
\end{align*}
Therefore $y_1 \notin Q$ (since $(y_1\biss y_r)$ is part of a
reducing system of parameters of $M$), a contradiction.

(ii) $\Rightarrow$ (i): Let $\ov M:=M/x_1M$ and take $P \in \supp
\ov M \cap V(x_2\biss x_r) = \supp M \cap V(x_1\biss x_r)$ with
$\dim R/P = \dim \ov M\, - (r-1) = \dim M\, -r$. Then $M_P$ is an
$r$-dimensional \cm\ module (over $R_P$) by our assumption and
$(x_1\biss x_r)$ is a system of parameters of $M_P$ and hence a
regular sequence on $M_P$ by Remark \ref{4}(1). But then $\ov M_P
\cong M_P/x_1M_P$ is an ($r-1$)-dimensional \cm\ module (over
$R_P$).

By the induction hypothesis $(x_2\biss x_r)$ is part of a reducing
system of parameters of $\ov M$ and therefore it remains to show
that $x_1 \notin Q$ for all $Q \in \ass M$ with $\dim R/Q = d-1$.

Suppose this is not the case. Choose $Q \in \ass M$ with $\dim R/Q
= d-1$ and $x_1 \in Q$. By Corollary \ref{cor} there is a prime $P
\in \supp M\cap V(x_2\biss x_r)$ such that $Q \subseteq P$ and
$\dim R/P \ge d-1 -(r-1) = d-r$. But then $P \in \supp M\cap
V(x_1\biss x_r)$ (since $x_1 \in Q \subseteq P$) and therefore
$\dim R/P \le d-r$, i.e. $\dim R/P = d-r$. By our assumption,
$M_P$ is an $r$-dimensional Cohen--Macaulay module. Since $QR_P
\in \ass M_P$ we therefore have
\begin{align*}
r &= \dim M_P = \dim R_P/QR_P = \dim (R/Q)_P\\[.3pc] &\le \dim R/Q -
\dim R/P = d-1 -(d-r)\\[.3pc] &= r-1
\end{align*}
(see Remark \ref{4}(2)), a contradiction.\hfill $\Box$
\end{proof}

\setcounter{coro}{14}
\begin{coro}\label{cor2}$\left.\right.$\vspace{.5pc}

\noindent Let $(x_1\biss x_r)$ be part of a reducing system of
parameters of $M$. If $r < d${\rm ,} then $(x_{\pi(1)}\biss
x_{\pi(r)})$ is part of a reducing system of parameters of $M$ for
any permutation $\pi$ of $\{1\biss r\}$.
\end{coro}

We note that the statement of this corollary is not true in
general if $r=d$, see the following Example \ref{ex}.

\setcounter{exa}{15}
\begin{exa}\label{ex}
{\rm Let $R := K\,[\mspace{-3mu}[ X,Y,Z]\mspace{-3mu}]$, where $K$
is a field and $X, Y, Z$ are indeterminates. For
\begin{equation*}
M := R/(XY,XZ)R\quad \mbox{ and }\quad x_1:=Y, \: x_2 := X+Y+Z,
\end{equation*}
$(x_1,x_2)$ is a system of parameters of $M$, but not a reducing
system of parameters. $(x_2,x_1)$ is a reducing system of
parameters of $M$ (not a regular sequence of $M$).}
\end{exa}

Finally we define the following.

\setcounter{definit}{16}
\begin{definit}\label{defcm}

{\rm \begin{align*} \mathcal{CM}\,(M) &:= \{P\in \supp M |\dim
R/P + \dim M_P = d \;\mbox{ and }\; M_P \mbox{ is}\\[.3pc]
&\qquad \mbox{a \cm\ module over } R_P \}
\end{align*}
(the {\it strong \cm\ locus of} $\supp M$) and for $0 \le r \le d$
\begin{equation*}
\mathcal{CM}_r(M) := \{P\in \mathcal{CM}\,(M)|\dim M_P = r\}.
\end{equation*}}
\end{definit}

\setcounter{rem}{17}
\begin{rem}$\left.\right.$

{\rm \begin{enumerate}
\renewcommand{\labelenumi}{(\arabic{enumi})}
\leftskip .15pc
\item We have
\leftskip .4pc \vspace{6pt}

(i) $\mathcal{CM}_0(M)=\assh M$ and, if $d \ge 1$,

(ii) $\mathcal{CM}_1(M) = \{P\in \supp M|\dim R/P = d-1\}\setminus
\ass M$,

(iii) $\mathcal{CM}\,(M) = \bigcup_{r=0}^d
\mathcal{CM}_r(M)$.\vspace{1pc}
\item The following conditions are equivalent
\leftskip .5pc \vspace{6pt}

(i) $M$ is a \cm\ module,

(ii) $\mathcal{CM}\,(M) = \supp M$,

(iii) $\m \in \mathcal{CM}\,(M)$.\vspace{1pc}
\item If $\supp M$ is equidimensional and catenarian then $\mathcal{CM}\,(M)$ coincides with the ordinary \cm\ locus of $\supp M$. This is the case, for example, when $\dim M \le 1$ or when $R$ is an epimorphic image of a local \cm\ ring and $M$ is
equidimensional.
\end{enumerate}}
\end{rem}

\setcounter{propo}{18}
\begin{propo}$\left.\right.$\vspace{.5pc}

\noindent For $r \in \mathbb N${\rm ,} $r < d${\rm ,} we have
\begin{align*}
\mathcal{CM}_r(M) &= \{P|P \in \ass M/(x_1\biss x_r)M,\; \dim R/P
= d-r,\\[.3pc]
&\quad\, (x_1\biss x_r) \mbox{ part of a reducing system of
parameters of } M\}.
\end{align*}
\end{propo}

\begin{proof}
By Theorem \ref{10} we have `$\supseteq$' and equality holds
(trivially) for $r=0$. Therefore it remains to verify the validity
of the inclusion `$\subseteq$' for $r\ge 1$.\pagebreak

Let $P \in \mathcal{CM}_r(M)$. Since $M_P$ is a \cm\ module with
$\dim M_P = r \ge 1$, we have $PR_P \notin \ass M_P$ and hence $P
\notin \ass M$. Moreover, $\dim R/P = d -\dim M_P = d-r$.

Let $Q \in \ass M$ with $\dim R/Q \ge d-1$. Then $P \not\subseteq
Q$ since $P=Q$ is impossible ($P \notin \ass M$) and $P \subset Q$
would imply $\dim R/P = d$ contradicting again `$P \notin \ass
M$'. Therefore we can find an $x_1 \in P$ with $x_1 \notin Q$ for
all $Q \in \ass M$ with $\dim R/Q \ge d-1$. By construction,
$(x_1)$ is part of a reducing system of parameters of $M$ and a
regular sequence on $M_P$ by Lemma \ref{5} and Remark \ref{4}(2).

If $r > 1$ we continue this procedure by passing to $M/x_1M$ and
we can construct elements $x_1\biss x_r\in P$ inductively on $r$
such that $(x_1\biss x_r)$ forms a part of a reducing system of
parameters of $M$ and a regular sequence on $M_P$. Let $\ov M :=
M/(x_1\biss x_r)M$. Since $\dim \ov M_P = \dim M_P/(x_1\biss
x_r)M_P = \dim M_P - r = 0$, $P$ is minimal in $\supp \ov M$ and
therefore $P \in \ass \ov M$.\hfill $\Box$
\end{proof}

\end{document}